\magnification=1200
\vsize=185mm
\hsize=135mm
\voffset=5mm
\null

\def\hypBasicAncA{1.1}
\def\hypBasicAncC{1.2}
\def\hypBasic{1.3}
\def \hypCommonSupp{1.4}
\def \HypThreeDd{1.5}
\def\Hypcomega{1.6}
\def \HypThreeDc{1.7}
\def \HypThreeDb{1.8}
\def\HypThreeDbb{1.9}

\def\ClaimWeak{2.1}
  \def\EstimElambdaA{2.2}
 \def\EstimElambdaB{2.3}
 \def\EstimElambdaBb{2.4}
  \def\EstimElambdaBc{2.5}
 \def\condGammaInterpol{2.6}
  \def\Deffu{2.7}
 \def\estimulocal{2.8} 
\def\estimGraduDelta{2.9}
\def\estimLpDelta{2.10}
\def\estimLoneMinus{2.11}
\def\LowerUHopf{2.12}
\def\estimfLoneDelta{2.13}
\def\estimDeltaGamma{2.14}
\def\estimcunhalfA{2.15}
\def\IneqSobolev{2.16}
\def\IneqHardy{2.17}
\def\estimGraddeltaomega{2.18}
\def\estimHardyomega{2.19}
\def\estimSobolevomega{2.20}
\def\estimNhalfomega{2.21}
\def\estimGradBB{2.22}
\def\estimHardyB{2.23}
\def\estimLoneMinusB{2.24}

\def\EqnUnbddA{3.1}
\def\EqnUnbddB{3.2}

\def\equivHopf{A.1}
\def\ControlNablaVarphi{A.2}
\def\ClaimAaa{A.3}
\def\ClaimAab{A.4}
\def\ClaimAa{A.5}
\def\ineqAa{A.6}
\def\ineqAb{A.7}
\def\ineqAbb{A.8}
\def\ineqAc{A.9}
\def\ClaimBa{A.10}

\def\ADP{1}
\def\AP{2}
\def\AL{3}
\def\ACJT{4}
\def\Av{5}
\def\BVV{6}
\def\BMPa{7}
\def\BMPb{8}
\def\BC{9}
\def\BT{10}
\def\CRT{11}
\def\FM{12}
\def\FPa{13}
\def\FPb{14}
\def\FSW{15}
\def\FMT{16}
\def\GM{17}
\def\JS{18}
\def\MPS{19}
\def\Ne{20}
\def\OpKu{21}
\def\QSa{22}
\def\QSb{23}

\def\eps{\varepsilon}

\def \trait (#1) (#2) (#3){\vrule width #1pt height #2pt depth #3pt}
\def \fin{\null\hfill
        \trait (0.1) (5) (0)
        \trait (5) (0.1) (0)
        \kern-5pt
        \trait (5) (5) (-4.9)
        \trait (0.1) (5) (0)
\medskip}
\null

\newtoks\hautpagegauche
\newtoks\hautpagedroite
\newtoks\paragraphecourant
\newtoks\chapitrecourant
\hautpagegauche={}
\hautpagedroite={}
\headline={\ifodd\pageno\the\hautpagedroite\else\the\hautpagegauche\fi}

\font\TenEns=msbm10
\font\SevenEns=msbm7
\font\FiveEns=msbm5
\newfam\Ensfam
\def\Ens{\fam\Ensfam\TenEns}
\textfont\Ensfam=\TenEns
\scriptfont\Ensfam=\SevenEns
\scriptscriptfont\Ensfam=\FiveEns
\def\R{{\Ens R}}

\def\eps{\varepsilon}


\font\itsmall=cmsl9

\font\eightrm=cmr9
\font\sixrm=cmr6
\font\fiverm=cmr5

\font\eighti=cmmi9
\font\sixi=cmmi6
\font\fivei=cmmi5

\font\eightsy=cmsy9
\font\sixsy=cmsy6
\font\fivesy=cmsy5

\font\eightit=cmti9
\font\eightsl=cmsl9
\font\eighttt=cmtt9

\def\eightpoint{\def\rm{\fam0\eightrm}
\textfont0=\eightrm
\scriptfont0=\sixrm
\scriptscriptfont0=\fiverm

\textfont1=\eighti
\scriptfont1=\sixi
\scriptscriptfont1=\fivei

\textfont2=\eightsy
\scriptfont2=\sixsy
\scriptscriptfont2=\fivesy

\textfont3=\tenex
\scriptfont3=\tenex
\scriptscriptfont3=\tenex

\textfont\itfam=\eightit \def\it{\fam\itfam\eightit}
\textfont\slfam=\eightsl \def\it{\fam\slfam\eightsl}
\textfont\ttfam=\eighttt \def\it{\fam\ttfam\eighttt}
}


\font\pc=cmcsc9
\font\itsmall=cmsl9

\def \trait (#1) (#2) (#3){\vrule width #1pt height #2pt depth #3pt}
\def \fin{\hfill
       \trait (0.1) (5) (0)
       \trait (5) (0.1) (0)
       \kern-5pt
       \trait (5) (5) (-4.9)
       \trait (0.1) (5) (0)
\medskip}
%


\paragraphecourant={\rmt ELLIPTIC EQUATION WITH CRITICAL GROWTH}
\chapitrecourant={\rmt Ph. SOUPLET}
\footline={\ifnum\folio=1 \hfill\folio\hfill\fi}
\hautpagegauche={\tenrm\folio\hfill\the\chapitrecourant\hfill}
\hautpagedroite={\ifnum\folio=1 \hfill\else\hfill\the\paragraphecourant\hfill\tenrm\folio\fi}

\font\rmb=cmbx8 scaled 1125 \rm

\centerline{\rmb A PRIORI ESTIMATES AND BIFURCATION OF SOLUTIONS}
\vskip 0.5mm
\centerline{\rmb FOR A NONCOERCIVE ELLIPTIC EQUATION}
\vskip 0.5mm
\centerline{\rmb  WITH CRITICAL GROWTH IN THE GRADIENT}

\vskip 5mm
\centerline{\pc Philippe SOUPLET}
\vskip 3mm
\centerline{\itsmall Universit\'e Paris 13, Sorbonne Paris Cit\'e, CNRS UMR 7539}
\centerline{\itsmall Laboratoire Analyse G\'eom\'etrie et Applications}
\centerline{\itsmall 93430 Villetaneuse, France. Email: souplet@math.univ-paris13.fr}

\vskip 4mm

\baselineskip=12pt
\font\rmt=cmr9

\setbox1=\vbox{
\hsize=120mm
{\baselineskip=11pt \parindent=3mm \eightpoint \rmt
{\pc Abstract:}\ 
 We study nonnegative solutions of the boundary value problem
 $$-\Delta u = \lambda c(x)u + \mu(x)|\nabla u|^2 + h(x),\quad
 u\in H^1_0(\Omega)\cap L^\infty(\Omega),
 \leqno(P_\lambda)$$
where $\Omega$ is a smooth bounded domain,
 $\mu, c\in L^\infty(\Omega)$, $h\in L^r(\Omega)$ for some $r > n/2$ and
$\mu,c,h > {\hskip -3.5mm} {\atop \neq} 0$.
Our main motivation is to study the ``noncoercive'' case.
Namely, unlike in previous work on the subject, we do not assume $\mu$ to be positive everywhere in $\Omega$.

 In space dimensions up to $n=5$, we establish uniform a priori estimates for weak solutions of ($P_\lambda$) when $\lambda>0$ is bounded away from $0$.
This is proved under the assumption that the supports of $\mu$ and $c$ intersect,
a condition that we show to be actually necessary,
and in some cases we further assume that $\mu$ is uniformly positive on the support of $c$
and/or some other conditions.

As a consequence of our a priori estimates, assuming that ($P_0$) has a solution, we deduce the existence of a continuum ${\cal C}$ of solutions,
such that the projection of ${\cal C}$ onto the $\lambda$-axis is an interval of the form $[0,a]$ for some $a>0$ and that the continuum ${\cal C}$ bifurcates from infinity to the right of the axis $\lambda=0$. In particular, for each $\lambda>0$ small enough, 
problem $(P_\lambda)$ has at least two distinct solutions.

\vskip 0.2cm

{\pc Keywords:}\
Elliptic equation, critical growth in the gradient, a~priori estimates, existence, multiplicity, bifurcation, $L^p_\delta$ spaces, weighted Sobolev inequalities, Hardy inequalities

}}

\hskip 2mm \hbox{\box1}

\bigskip

{\bf 1. Introduction and main results.}
\medskip

In this article, we consider the following Dirichlet problem:
$$\left\{\eqalign{
-\Delta u&=\mu(x)|\nabla u|^2+\lambda{\hskip 1pt} c(x)u+h(x),\cr
\noalign{\vskip 1mm}
 u&\in H^1_0(\Omega)\cap L^\infty(\Omega).
}\right.
\leqno(\hbox{$P_\lambda$})$$
Here $\Omega\subset \R^n$ is a bounded domain of class $C^2$,
$\mu, c, h$ are given functions, whose regularity will be specified below, and $\lambda$ is a real parameter.
By a solution, we mean a weak solution in the sense of the usual integral formulation,
with test-functions in $H^1_0(\Omega)\cap L^\infty(\Omega)$.
It is known (see [\ACJT]) that, under assumption (\hypBasic) below, any solution $u$ of $(P_\lambda)$
is H\"older continuous in~$\overline\Omega$.

\smallskip

Elliptic equations with a gradient dependence up to the critical (quadratic) growth were studied by Boccardo, Murat and Puel in the 80's and a large literature on the subject has appeared since then.
Many results are known in the case $\lambda=0$ 
(see e.g. [\AL, \MPS, \AP, \FPa, \GM,  \FPb, \FM]) or $\lambda<0$ (see e.g. [\BMPa, \BMPb]).
We shall here consider the case $\lambda>0$ and will be concerned with questions of a priori estimates, existence, multiplicity and bifurcation of solutions.
\smallskip

Denote the nonnegative solution set 
$$\Sigma=\bigl\{(\lambda,u)\in [0,\infty)\times C(\overline\Omega);\  \hbox{ $u\ge 0$ and $u$ solves } (P_\lambda)\bigr\}.$$
Interesting properties of the set $\Sigma$ were recently established in [\ACJT] under the assumptions (with $n\ge 3$):
$$\leqalignno{
&\hbox{$\mu \in L^\infty(\Omega),\quad c, h\in L^r(\Omega)$\ for some $r>n/2$,
\ \ $c, h > {\hskip -3.5mm} {\atop \neq} 0$,} 
&(\hypBasicAncA)\cr
&\hbox{$\mu(x)\ge \mu_0>0$.}&(\hypBasicAncC)
}
$$
For every $\eps>0$, it is shown that nonnegative solutions of $(P_\lambda)$ with $\lambda\ge \eps$ satisfy a uniform a priori estimate. 
Next assume in addition that ($P_0$) has a (necessarily nonnegative) solution
(see Remark~1.1 below for known sufficient conditions).
Then it is shown in [\ACJT] that there exists a continuum ${\cal C}\subset \Sigma$,
such that the projection of ${\cal C}$ onto the $\lambda$-axis is an interval of the form $[0,a]$ for some $a>0$ and that the continuum ${\cal C}$ bifurcates from infinity to the right of the axis $\lambda=0$. In particular, for each $\lambda>0$ small enough, ${\cal C}$ contains at least two distinct solutions of $(P_\lambda)$.
\smallskip

We note that, when $\mu(x)=\mu$ is a positive constant, multiplicity results 
(actually up to an explicit value $\lambda_0>0$ of $\lambda$) have been obtained before in [\JS], using the transformation $v=e^{\mu u}$ 
(see also [\ADP]).
We stress that the multiplicity results 
in [\ACJT] allow nonconstant functions $\mu(x)$, in which case such a transformation is not available.
However, the coercivity of the function $\mu(x)$, i.e.~(\hypBasicAncC), is still needed in [\ACJT].
\smallskip

Our main goal here is to establish similar results as in [\ACJT] for
noncoercive functions $\mu(x)$,
namely to allow just $\mu\ge 0$, $\mu\not\equiv 0$
(possibly at the expense of additional assumptions on~$c$).
The main difficulty here is to establish a priori estimates without assuming (\hypBasicAncC).
Indeed, this assumption seems necessary in [\ACJT], in order to apply the method 
of Brezis and Turner~[\BT] based on Hardy-Sobolev inequalities.
Therefore, we need some new ideas
(see after the statement of the results, 
for a brief description of the main arguments of our proofs).

\medskip

We assume:
$$\mu, c \in L^\infty(\Omega),\quad h\in L^r(\Omega)\ \hbox{for some $r>\max(1,n/2)$,} 
\quad\hbox{$ \mu, c, h > {\hskip -3.4mm} {\atop \neq} 0.$}
\leqno(\hypBasic)$$
Also, we shall  make the following essential assumption
of {\it intersecting supports} for $\mu$ and~$c$:
$$\mu, c\ge \eta \hbox{ on $B(x_0,\rho)\subset\Omega$ for some $\rho, \eta>0.$}
\leqno(\hypCommonSupp)$$
In low dimensions $n\le 2$, this turns out to be sufficient to guarantee a priori estimates.

\proclaim Theorem 1. 
Let $n\le 2$ and assume (\hypBasic), (\hypCommonSupp).
Then for any $\Lambda_1 > 0$ there exists a constant $M > 0$ such that, for each 
$\lambda\ge\Lambda_1$, any nonnegative solution of ($P_\lambda$) satisfies $\|u\|_\infty\le M$.

We shall see right away that the assumption (\hypBasic) of intersecting supports for $\mu$ and $c$ is also essentially necessary.

\proclaim Theorem 2.
Consider problem ($P_\lambda$) with 
$$n=1,\quad \Omega=(0,3),\quad \mu(x)=\chi_{(1,2)},\quad
c(x)=\chi_{(2,3)},\quad h=0.$$
 There exists a sequence $\lambda=\lambda_j\to \pi^2/4$
and a sequence $u_j\in H^2(\Omega)\cap H^1_0(\Omega)$ of solutions of ($P_\lambda$)
such that $\|u_j\|_\infty\to\infty$ as $j\to\infty$.

We now turn to the higher dimensional range $3\le n\le 5$.
In this case, beside (\hypCommonSupp), we need additional assumptions
to guarantee a priori estimates.
In our next result, we shall assume, roughly speaking,  that $\mu$ 
{\it is positively bounded below on the support of } $c$. 
However, we do not know presently whether this assumption is technical or not
(nor the assumptions in Theorem 4 below).
In particular, we do not know if the a priori estimates are still
true in dimensions $n\ge 6$. 
The restriction $n\le 5$ comes from the fact that our method (see at the end
of this section)
requires to estimate the function $c(x)u$ in $L^p$ for some $p>n/2$, whereas
the $L^p$ estimates that we are able to derive are limited to $p\le 2n/(n-2)$,
due to Sobolev imbeddings or to even more stringent functional inequalities
(note that $2n/(n-2)$ and $n/2$ precisely coincide for $n=6$).

\proclaim Theorem 3.
Let $3\le n\le 5$ and let (\hypBasic), (\hypCommonSupp) be satisfied.
Assume that there exists a $C^2$ domain $\omega\subset\Omega$ and a constant $\mu_0>0$, such that
$$
\hbox{$\mu\ge\mu_0$ on $\omega$
and
${\rm Supp}(c)\subset\overline\omega$.}
\leqno(\HypThreeDd)$$
If $n=5$, assume in addition that
$$c(x)\le C_1[{\rm dist}(x,\partial\omega)]^\sigma,\ \ x\in \omega,
\quad\hbox{ for some $C_1,\sigma>0$.}
\leqno(\Hypcomega)$$
Then for any $\Lambda_1 > 0$ there exists a constant $M > 0$ such that, for each 
$\lambda\ge\Lambda_1$, any nonnegative solution of ($P_\lambda$) satisfies $\|u\|_\infty\le M$.

\medskip
As usual, ${\rm Supp}$ is here understood in the sense of essential support.
As a special case of Theorem~3, we see that the a priori estimate holds for instance if $3\le n\le 5$, $c$ is compactly supported and $\mu\ge\mu_0>0$ on a neighborhood of $c$.
\medskip

We now give additional results in the case of dimension $n=3$, which is rather special.
Indeed, without assuming (\HypThreeDd),
we can then obtain a priori estimates under various, relatively mild, assumptions
either on $\mu$ or $c$.

\proclaim Theorem 4.
Let $n=3$ and (\hypBasic), (\hypCommonSupp) be satisfied.
Assume in addition that either
$$c(x)\le C_1[{\rm dist}(x,\partial\Omega)]^\sigma,\ \ x\in \Omega,
\quad\hbox{ for some $C_1,\sigma>0$,}
\leqno(\HypThreeDc)$$
or
$${\rm Supp}(\mu)\subset\subset\Omega,
\leqno(\HypThreeDb)$$
or
$$\mu(x)\ge c_1[{\rm dist}(x,\partial\Omega)]^\sigma \quad\hbox{ on a neighborhood of $\partial\Omega$ for some $c_1>0$ and $\sigma<2$.}
\leqno(\HypThreeDbb)$$
Then for any $\Lambda_1 > 0$ there exists a constant $M > 0$ such that, for each 
$\lambda\ge\Lambda_1$, any nonnegative solution of ($P_\lambda$) satisfies $\|u\|_\infty\le M$.

\medskip

As a consequence of Theorems 1, 3, 4 and of [\ACJT, Theorem 1.2]
(see also the proof of Theorem 1.3 in [\ACJT]),
one deduces the following result on existence, multiplicity and bifurcation.

\proclaim Theorem 5. Assume that ($P_0$) has a solution
and let the assumptions of Theorem~1 (resp. Theorem~3, 4) be in force.
Then there exists a continuum ${\cal C}\subset \Sigma$,
such that the projection of ${\cal C}$ onto the $\lambda$-axis is an interval of the form $[0,a]$ for some $a>0$ and that the continuum ${\cal C}$ bifurcates from infinity to the right of the axis $\lambda=0$. In particular, for each $\lambda>0$ small enough, 
problem $(P_\lambda)$ has at least two distinct nonnegative solutions.

\medskip

\medskip

{\bf Remarks 1.1.} (a) The existence of a solution for ($P_0$)
is known under various smallness assumptions. For instance, in [\FM], it is assumed that
$$\|\mu\|_\infty\|h\|_{n/2}<S_n^2
=\inf\left\{{\|\nabla \phi\|_2\over \|\phi\|_{2^*}}\,; 0\not\equiv\phi\in H^1_0(\Omega)\right\}$$
(where $S_n$ is the best constant in the Sobolev inequality and $n\ge 3$).
In the particular case $\mu(x)= \mu_0 > 0$ and $h\ge 0$,
a more precise sufficient condition for the existence of a solution for ($P_0$) is given [\ADP] by
$$\mu<\inf\left\{{\int_\Omega |\nabla \phi|^2\,dx\over \int_\Omega h(x)\phi^2\,dx}\,; 
0\not\equiv\phi\in H^1_0(\Omega)\right\}.$$
\smallskip

(b) We stress that some of our results could be extended to the case
when the function $c$ belongs to $L^p$ for suitable $p$ (instead of $L^\infty$).
However, since our main motivation here is to treat the case of noncoercive $\mu$,
 we have left this aside from simplicity.
 
 \medskip
 
 To prove a priori estimates without the coercivity assumption (\hypBasicAncC), we shall combine various ingredients.
We first reduce the desired uniform bounds to estimating the term $c(x)u$ in a suitable $L^q$ space (namely, $q>n/2$; see Proposition~1). This is done by using exponential test-functions and auxiliary unknowns, in the spirit of, e.g., 
Boccardo-Murat-Puel [\BMPa].
To derive the necessary $L^q$ estimates, we first get (Lemma~5)
some basic weighted 
$L^1$ estimates for $u$ and for the right-hand side of $(P_\lambda)$.
This is obtained by using the quantitative Hopf Lemma 
in Brezis-Cabr\'e [\BC], and a test-function which solves a singular auxiliary problem (originally from Crandall-Rabinowitz-Tartar [\CRT]). We then apply smoothing effects in $L^p_\delta$, the Lebesgue spaces weighted by the distance to the boundary,
along with suitable interpolation.
Note that such smoothing effects were used before for problems without gradient terms
(see [\BVV, \FSW, \QSa]) but not, as far as we know, for problems with (critical)
gradient terms.
 In higher dimensional cases, we also need to take advantage 
 of the weighted $H^1$ bound guaranteed by the estimate of the right-hand side.
 To this end, we rely on suitable weighted Sobolev and Hardy inequalities
(different from those in [\ACJT] or [\BT]); see Proposition~2 in section 2.4 below.

\medskip
{\bf Acknowledgement.} The author thanks Louis Jeanjean for stimulating discussion
concerning this work.

\medskip

{\bf 2. Proofs of a priori estimates.}
\medskip

In the rest of the article, we denote by 
$$\delta(x)=\delta_\Omega(x)={\rm dist}(x,\partial\Omega)$$
 the function distance to the boundary
(the subscript $\Omega$ will be dropped when no confusion arises).
Conjugate exponents will be denoted by ${}'$ (i.e., $p'=p/(p-1)$, $1<p<\infty$).
For $1\le p\le\infty$, the $L^p(\Omega)$ and $L^p_\delta(\Omega)$ norms will be denoted by
$\|\cdot\|_p$ and $\|\cdot\|_{p,\delta}$.
Recall that
$$\|u\|_{p,\delta}=\Bigr(\int_\Omega |u(x)|^p\delta(x)\, dx\Bigl)^{1/p},\quad 1\le p<\infty.$$
The notation will be also used for $0<p<1$ although it is not a norm in this case.

Moereover, we recall (see [\ACJT]) that problem ($P_\lambda$) admits no 
nontrivial nonnegative solutions
for $\lambda> \gamma_1$, where $\gamma_1>0$ is the first eigenfunction of the
problem 
$$-\Delta \varphi=\gamma_1c(x)\varphi,\quad\varphi\in H^1_0(\Omega).$$

\medskip
\ \ {\bf 2.1. Reduction to a suitable $L^q$ estimate.}
\medskip

The goal of this subsection is to reduce the proof of $L^\infty$ estimate
to a suitable $L^q$ estimate.
We shall prove:

\proclaim Proposition 1. Assume (\hypBasic). 
Set $q=n/2$ if $n\ge 3$ or fix any $q\in (1,\infty)$ if $n\le 2$.
Let $\eps,\alpha,M>0$. There exists $C_1>0$ such that, 
for any $\lambda\in [0,\gamma_1]$ and any nonnegative solution of ($P_\lambda$),
$$\|cu^{1+\eps}\|_{q}+\|u\|_\alpha\le M
\Longrightarrow
\|u\|_\infty\le C_1.$$

\medskip
The proposition is a direct consequence of the following two lemmas.
\medskip

\proclaim Lemma 1. Assume (\hypBasic). There exists $K=K(r,n,\|\mu\|_\infty)>0$
with the following property. For any $M_1>0$,
there exists $C_1>0$ such that,
for any $\lambda\in [0,\gamma_1]$ and any nonnegative solution of ($P_\lambda$),
$$\|e^{Ku}\|_1\le M_1
\Longrightarrow
\|u\|_\infty\le C_1.$$

\medskip

\proclaim Lemma 2. Assume (\hypBasic). Let $q$ be as in Proposition 1
and let $\eps,\alpha,k,M>0$. There exists $M_1>0$ such that,
for any $\lambda\in [0,\gamma_1]$ and any nonnegative solution of ($P_\lambda$),
$$\|cu^{1+\eps}\|_{q}+\|u\|_\alpha\le M
\Longrightarrow
\|e^{ku}\|_1\le M_1.$$

\medskip

{\it Proof of Lemma 1.} 
For any $k>0$, we have $e^{k u}-1\in H^1_0(\Omega)$,
due to $u\in H^1_0(\Omega)\cap L^\infty(\Omega)$.
We claim that
$$-\Delta (e^{k u}-1)= k e^{k u}\bigl((\mu(x)-k)|\nabla u|^2+\lambda{\hskip 1pt} c(x)u+h(x)\bigr)
\ \ \hbox{in the weak sense in $\Omega$}
\leqno(\ClaimWeak)$$
(i.e., with test functions in 
$H^1_0(\Omega)\cap L^\infty(\Omega)$).
Indeed, setting $w=e^{k u}-1$ and 
$F=\mu |\nabla u|^2+\lambda{\hskip 1pt} cu+h\in L^1(\Omega)$,
we compute, for any $\phi\in H^1_0(\Omega)\cap L^\infty(\Omega)$,
$$\int_\Omega \nabla w\cdot\nabla\phi
=\int_\Omega ke^{ku}\nabla u\cdot\nabla\phi
=\int_\Omega \Bigl(\nabla (ke^{ku}\phi)-k^2e^{ku}\phi\nabla u\Bigr)\cdot\nabla u.$$
On the other hand, using $ke^{ku}\phi\in H^1_0(\Omega)\cap L^\infty(\Omega)$ as a test-function in 
($P_\lambda$), we have
$$\int_\Omega\nabla (ke^{ku}\phi)\cdot\nabla u=\int_\Omega ke^{ku} \phi F.$$
It follows that
$$\int_\Omega \nabla w\cdot\nabla\phi
=\int_\Omega ke^{ku}\bigl(F-k|\nabla u|^2\bigr)\phi,$$
hence the claim (\ClaimWeak).

Now choosing $k=\|\mu\|_\infty$, as a consequence of (\ClaimWeak), we have
$$-\Delta (e^{k u}-1) \le g:=k e^{k u}\bigl(\lambda{\hskip 1pt} c(x)u+h(x)\bigr)$$
in the weak sense in $\Omega$.
Since $r>\max(1,n/2)$, we may fix $r_1\in (\max(1,n/2),r)$ and $r_2\in (1,\infty)$ such that
$\textstyle{1\over r_1}=\textstyle{1\over r}+\textstyle{1\over r_2}$.
Using the maximum principle, the standard $L^{r_1}$-$L^\infty$ estimate for the Dirichlet Laplacian with $r_1>n/2$ (see, e.g., [\QSb, Proposition 47.5]), H\"older's inequality and $e^u\ge u$, we then obtain
$$\|e^{k u}-1\|_\infty\le C(r_1,\Omega) \|g\|_{r_1}\le 
kC(r_1,\Omega)\bigl[\lambda\|c\|_\infty\|e^{(k+1)u}\|_{r_1}+\|e^{ku}\|_{r_2}\|h\|_r\bigr].$$
The Lemma follows with $K=\max((k+1)r_1,kr_2)$. \fin
 
 \medskip

{\it Proof of Lemma 2.} If $n\le 2$, we may assume $q<r$ without loss of generality.
Fix $M>0$ and assume
$$\|cu^{1+\eps}\|_{q}+\|u\|_\alpha\le M.
\leqno(\EstimElambdaA)$$
For $k>0$, testing problem ($P_\lambda$) with $e^{k u}-1
\in H^1_0(\Omega)\cap L^\infty(\Omega)$, we get
$$\int_\Omega k e^{k u}|\nabla u|^2=\int_\Omega (e^{k u}-1)
\bigl(\mu(x)|\nabla u|^2+\lambda{\hskip 1pt} c(x)u+h(x)\bigr).$$
Assume $k\ge 1+\|\mu\|_\infty$ without loss of generality. Then
$$\int_\Omega  e^{k u}|\nabla u|^2\le\int_\Omega (e^{k u}-1)(\lambda{\hskip 1pt} c(x)u+h(x)).$$
Set $s=q'$ and note that $s=n/(n-2)$ if $n\ge 3$.
By Sobolev's inequality, we get
$$
\|e^{k u/2}-1\|_{2s}^2
\le C\int_\Omega  |\nabla(e^{k u/2}-1)|^2\le
C{k^2\over 4}\int_\Omega (e^{k u}-1)(\lambda{\hskip 1pt} c(x)u+h(x)).
\leqno(\EstimElambdaB)$$
Here and in the rest of the proof, $C$ denotes a generic constant independent of $u$, $\lambda$ and~$k$.
Also, using $X-1\le 2(\sqrt X-1)^2+C$ for all $X\ge 1$, we have
$$\|e^{k u}-1\|_{s}\le C\|(e^{k u/2}-1)^2\|_{s}+C
\le C\|e^{k u/2}-1\|_{2s}^2+C.
\leqno(\EstimElambdaBb)$$

Next, for any $\tau\in (0,1]$, by (\EstimElambdaB), (\EstimElambdaBb), $\lambda\le \gamma_1$ and $s'=q$, we have
$$\|e^{k u}-1\|_{s}-C
\le Ck^2\Bigl\|{e^{k u}-1\over u^\tau}\Bigr\|_{s} \bigl(\|cu^{1+\tau}\|_q+\|hu^{\tau}\|_q\bigr).\leqno(\EstimElambdaBc)$$
Recalling $q<r$ and choosing $\tau=\min(1,\eps,\alpha(r-q)/(rq))$,
we deduce from (\EstimElambdaA) that
$$\|hu^{\tau}\|_{q}^{q}\le \int_\Omega h^r+\int_\Omega u^{\tau rq/(r-q)}\le 
\|h\|_r^r+|\Omega|+\int_\Omega u^\alpha\le C+M^\alpha$$
and
$$\|cu^{1+\tau}\|_{q}\le \|c\|_{q}+\|cu^{1+\eps}\|_{q}\le C+M.$$
It thus follows from (\EstimElambdaBc) that
$$\|e^{k u}-1\|_{s}
\le k^2C_0(M)\Bigl\|{e^{k u}-1\over u^\tau}\Bigr\|_{s}+C.$$
where $C_0(M)$ is a positive constant depending on $M$
but otherwise independent of $u$, $\lambda$.

Now observe that, for any $A>0$,
$$\Bigl\|{e^{k u}-1\over u^\tau}\Bigr\|_{s}
\le \Bigl\|{e^{k u}-1\over u^\tau}\chi_{\{u^\tau\ge A\}} \Bigr\|_{s}
+\Bigl\|{e^{k u}-1\over u^\tau}\chi_{\{u^\tau< A\}} \Bigr\|_{s}
\le A^{-1} \bigl\|e^{k u}-1\bigr\|_{s} + C(A,k,\tau).$$
Choosing $A=2k^2C_0(M)$, we get
$$\|e^{k u}-1\|_{s}\le C_1(M,k,\alpha,\eps),$$
which implies the desired estimate.
\fin

For further reference, we give the following elementary interpolation Lemma.

\proclaim Lemma 3.
Let $\omega$ be any open subset of $\R^n$
and $\phi\in L^\infty(\omega)$ be such that $\phi>0$ a.e.
Let $1\le q\le m \le r<\infty$ and let $b,d\ge 0$
and $\gamma\in \R$ be such that
$$\gamma\ge d-{(r-m)(b+d)\over r-q},
\leqno(\condGammaInterpol)$$
We have, for any measurable function $v$ on $\omega$,
$$\int_\omega \phi^\gamma |v|^m\, dx
\le C\Bigl(\int_\omega \phi^{-b} |v|^q\, dx\Bigr)^\theta
\Bigl(\int_\omega \phi^{d} |v|^r\, dx\Bigr)^{1-\theta},
\qquad\theta={r-m\over r-q}
$$
(assuming both integrals on the RHS to be finite).

{\it Proof.} It suffices to notice that
$m=\theta q+(1-\theta)r$,
$\gamma\ge -\theta b+(1-\theta)d$
and to use H\"older's inequality and $\phi\in L^\infty(\omega)$.
\fin

\medskip
\ \ {\bf 2.2. Weighted $L^q$ estimates.}
\medskip

In the rest of section~2, 
when $u$ is a  solution of ($P_\lambda$)
we set
$$f=f_u:=\mu(x)|\nabla u|^2+\lambda{\hskip 1pt} c(x)u+h(x).
\leqno(\Deffu)$$
Also, $C$ will denote a generic constant independent of $u$ and of $\lambda\in [\Lambda_1,\gamma_1]$,
but possibly depending on all the other data and parameters
(the dependence on some particular parameters will be emphasized
if necessary).

We start with an integral a priori estimate of $u$ on a ball where both $\mu$ and $c$ are positively bounded below.

\proclaim Lemma 4. 
Assume (\hypBasic), (\hypCommonSupp) and $0<\Lambda_1\le\lambda\le\gamma_1$.
Then any nonnegative solution $u$ of ($P_\lambda$) satisfies
$$\int_{B_{\rho/2}(x_0)} e^{\eta u} \le C(\eta,\rho).
\leqno(\estimulocal)$$

{\it Proof.} Arguing similarly as in the proof of Lemma 1, we have, in the distribution sense,
$$-\Delta (e^{\eta u})= \eta e^{\eta u}\bigl((\mu(x)-\eta)|\nabla u|^2+\lambda{\hskip 1pt} c(x)u+h(x)\bigr)
\ge \Lambda_1\eta^2 e^{\eta u}u\quad\hbox{in $B:=B_{\rho}(x_0)$},$$
where we used (\hypBasic) for the last inequality.
Therefore, for all $A>0$, there exists $C(A)>0$ such that
$$-\Delta (e^{\eta u})\ge Ae^{\eta u}-C(A)\quad\hbox{in $B$}.$$
Denote respectively by $\lambda_{1,B}$ and $\varphi_{1,B}$ the first Dirichlet eigenvalue and eigenfunction of $B$
and take $A=1+\lambda_{1,B}$.
Testing this inequality with $\varphi_{1,B}\in H^1_0(B)\cap C(\overline B)$,
we easily obtain 
$$\int_{B_{\rho/2}(x_0)} e^{\eta u} \le C \int_{B} e^{\eta u}\varphi_{1,B}\le C.$$
\fin

The next lemma gives our primary a priori estimates for $u$.

\proclaim Lemma 5. 
Assume (\hypBasic), (\hypCommonSupp)
and $0<\Lambda_1\le\lambda\le\gamma_1$.
Then any nonnegative solution $u$ of ($P_\lambda$) satisfies
the following a priori estimates:
$$\int_\Omega |\nabla u|^2\mu(x)\delta(x)\,dx\le C,
\leqno(\estimGraduDelta)$$
$$\|u\|_{p,\delta}\le C(p),\quad p<(n+1)/(n-1),
\leqno(\estimLpDelta)$$
$$\|\delta^{-\gamma}u\|_{1}\le C(\gamma),\quad 0<\gamma<1.
\leqno(\estimLoneMinus)$$

{\it Proof.} 
It is well known (see [\BC]) that
$$u\ge c_1\cdot\Bigl(\int_\Omega f\delta\Bigr) \delta\quad\hbox{in $\Omega$}.
\leqno(\LowerUHopf)$$
for some constant $c_1=c_1(\Omega)>0$.
From (\estimulocal), (\LowerUHopf), we obtain 
$$\|f\|_{1,\delta} = \int_\Omega f\delta\le C \Bigl( \int_{B_{\rho/2}(x_0)}u \Bigr) \Bigl(\int_{B_{\rho/2}(x_0)}\delta\Bigr)^{-1}\le C,
\leqno(\estimfLoneDelta)$$
hence in particular (\estimGraduDelta).
By the $L^1_\delta$-$L^p_\delta$ estimate for the Dirichlet Laplacian (see [\FSW]),
we deduce (\estimLpDelta) from (\estimfLoneDelta).

Next, for each $\gamma\in (0,1)$, it is known (see [\QSb, Lemma 10.4] and cf.~the references in [\QSb]), that there exists a function $\xi\in H^1_0(\Omega)$ 
such that $-\Delta \xi=\delta^{-\gamma}$ in ${\cal D}'(\Omega)$ and such that $\xi\le C\delta$ in $\Omega$. 
(Note that the RHS $\delta^{-\gamma}$ belongs to $L^1(\Omega)$; see after (\ClaimAa).)
Actually, for any $v\in H^1_0(\Omega)\cap L^\infty(\Omega)$, we have
$$\int_\Omega\nabla v\cdot\nabla\xi = \int_\Omega v\delta^{-\gamma}.$$
Testing ($P_\lambda$) with $\xi$ and using (\estimfLoneDelta), we obtain
$$\int_\Omega u\delta^{-\gamma}
=\int_\Omega\nabla u\cdot\nabla\xi =\int_\Omega f\xi\le C\int_\Omega f\delta\le C,$$
hence (\estimLoneMinus).

\fin

\ \ {\bf 2.3. Proof of Theorem 1 and of Theorem 4 under assumption
(\HypThreeDb) or~Ê(\HypThreeDc).}
\medskip

{\it Proof of Theorem 3 under assumption (\HypThreeDb).}
Since ${\rm Supp}(\mu)\subset\subset\Omega$ by assumption,
 estimate (\estimGraduDelta) implies
$\int_\Omega |\nabla u|^2\mu(x)\,dx\le C$. This combined with (\estimLoneMinus)
guarantees that $\|f\|_1 \le C$. 
We then deduce from the standard $L^1$-$L^p$ estimate for the Dirichlet Laplacian that
$$\|u\|_p\le C(p),\quad p<n/(n-2).$$
Since $n/(n-2)>n/2$ for $n=3$ and since $c\in L^\infty(\Omega)$,
the desired conclusion follows from Proposition 1.
\fin

{\it Proof of Theorem 1 and of Theorem 3 under assumption
(\HypThreeDc).}
We assume $n=2$ or $3$ (the case $n=1$ can be handled with 
obvious modifications).
Take $\gamma=0$ if $n=2$ and any $\gamma>0$ if $n=3$.
We claim that, for $\eps>0$ small,
$$\int_\Omega \delta^\gamma u^{{n\over 2}+\eps}\le C(\eps).
\leqno(\estimDeltaGamma)$$
To this end, we interpolate between (\estimLpDelta) and (\estimLoneMinus), applying Lemma 3 with $\phi=\delta_\Omega$,
$$m={n\over 2}+\eps,\quad q=1,\quad r={n+1\over n-1}-\eps,\quad b=1-\eps,\quad d=1,$$
with $\eps>0$ small
(note that $q\le m \le r$). Taking $\eps>0$ small enough, condition (\condGammaInterpol) is true provided
$$\gamma>1-2\Bigl({n+1\over n-1}-{n\over 2}\Bigr){n-1\over 2}
=1-{2(n+1)-n(n-1)\over 2}={n(n-3)\over 2}$$
and the claim follows.

We now easily deduce from (\estimDeltaGamma) that
$$\int_\Omega c^{n/2}u^{{n\over 2}+\eps}\le C.
\leqno(\estimcunhalfA)$$
Indeed, if $n=2$, this is true due to $\gamma=0$
and $c\in L^\infty(\Omega)$.
If $n=3$, since $c\le C_1\delta^\sigma$ by assumption (\HypThreeDc), 
this follows from (\estimDeltaGamma) with $\gamma=\sigma n/2>0$.

The desired conclusion is now a direct consequence of Proposition 1.
\fin

\ \ {\bf 2.4. Weighted Sobolev and Hardy inequalities
and proof of Theorem 3 and of Theorem 4 under assumption (\HypThreeDbb).}
\medskip

To go further, 
we want to exploit the weighted $H^1$ nature of estimate (\estimGraduDelta).
To this end, a key role will be played by the following
weighted Sobolev and Hardy inequalities.

\proclaim Proposition 2. 
Let $n\ge 2$ and let $\omega\subset \R^n$ be a bounded open set of class $C^2$
and denote $\delta(x)={\rm dist}(x,\partial\omega)$.
Let $1\le p<n$, $a>p-1$
and $k\ge0$.
\smallskip
(i) Set $p^*=np/(n-p)$. We have
$$\Bigl(\int_\omega \delta^{na/(n-p)}|v|^{p^*}\Bigr)^{p/p^*}
\le C\Bigl(\int_\omega \delta^k |v|\Bigr)^{p}
+C\int_\omega \delta^a |\nabla v|^p
\leqno(\IneqSobolev)$$
for all $v\in W^{1,p}(\omega)$.
\smallskip
(ii) We have
$$\int_\omega \delta^{a-p}|v|^{p}
\le C\Bigl(\int_\omega \delta^k |v|\Bigr)^{p}
+C\int_\omega \delta^a |\nabla v|^p
\leqno(\IneqHardy)$$
for all $v\in W^{1,p}(\omega)$.

Related results of Hardy and/or weighted Sobolev type 
have been known for a long time,
see for instance [\Ne, \Av, \OpKu, \FMT] and the references therein.
However, we haven't found a suitable reference for this specific statement.
Actually, the available results seem to involve $|v|^p$ instead of $|v|$ on the RHS,
which is not sufficient for our needs, or they may impose other restrictions
such as lowers bound on $k$
(although they can be more general in terms of weights or of domain regularity).
We thus give a proof in appendix.

Note that, since $u$ is not assumed to vanish on the boundary,
the restriction $a>p-1$ is necessary in assertion (ii)
(notice that if one would take $0\le a\le p-1$ with, for instance, $v\equiv 1$,
then the LHS integral would become divergent, 
contradicting the finiteness of the RHS).

\medskip

{\it Proof of Theorem 3.} 
Let $\omega$ be as in the statement of the theorem.
By (\HypThreeDd) and (\estimGraduDelta), we have
$$\int_\omega \delta_\omega |\nabla u|^2\le 
\int_\omega \delta_\Omega|\nabla u|^2\le C.
\leqno(\estimGraddeltaomega)$$
Applying Propostion 2(ii) with $p=2$, $a=1+\eps$, $k=0$, if follows from (\estimGraddeltaomega) and (\estimLoneMinus) that
$$\int_\omega \delta_\omega^{-1+\eps}u^2 \le C(\eps).
\leqno(\estimHardyomega)$$
If $n=3$, we directly deduce that, for $\eps>0$ small,
$$\int_\Omega c^{n/2}u^{{n\over 2}+\eps}\le 
C\int_\omega u^{{n\over 2}+\eps}\le C.$$

Next assume $n=4$ or $5$. Applying Proposition~2(i) with $p=2$, $a=1+\eps$, $k=0$, it follows
from (\estimGraddeltaomega) and (\estimLoneMinus) that
$$\int_\omega \delta_\omega^{n(1+\eps)/(n-2)}u^{2n/(n-2)} \le C.
\leqno(\estimSobolevomega)$$
We now interpolate between (\estimHardyomega) and (\estimSobolevomega).
To this end, we use Lemma 3 with $\phi=\delta_\omega$,
$$m={n\over 2}+\eps,\quad q=2,\quad r={2n\over n-2},\quad b=1-\eps,\quad 
d={n(1+\eps)\over n-2},$$
for $\eps>0$ small (note that $q\le m \le r$).  It follows that 
$$\int_\omega \delta_\omega^\gamma u^{{n\over 2}+\eps}\le C
\leqno(\estimNhalfomega)$$
whenever $\gamma$ satisfies condition (\condGammaInterpol).
Taking $\eps>0$ small enough, this is true provided
$$\gamma>
{n\over n-2}-{2(n-1)\over n-2}\Bigl({2n\over n-2}-{n\over 2}\Bigr){n-2\over 4}
={n\over n-2}-{n(6-n)\over n-2}{n-1\over 4}$$
i.e., $\gamma>n(n-5)/4$.

If $n=4$, we may take $\gamma=0$ in (\estimNhalfomega).
Choosing $\eps>0$ sufficiently small and
using ${\rm supp}(c)\subset \overline\omega$
and $c\in L^\infty(\Omega)$, we get
$$\int_\Omega c^{n/2}u^{{n\over 2}+\eps}\le 
C\int_\omega u^{{n\over 2}+\eps}\le C.$$

Finally, if $n=5$, by assumption (\Hypcomega), we have $c\le C_1\delta_\omega^\sigma\chi_{\overline\omega}$.
We may take $\gamma=\sigma n/2>0$ in (\estimNhalfomega).
Choosing $\eps>0$ sufficiently small, we get
$$\int_\Omega c^{n/2}u^{{n\over 2}+\eps}\le 
C\int_\omega \delta_\omega^{\sigma n/2} u^{{n\over 2}+\eps}\le C.$$
The conclusion follows from Proposition 1.
\fin

\medskip

\medskip

{\it Proof of Theorem 4 under assumption (\HypThreeDbb).}
By assumption, there exists $\eta>0$ such that the tubular neighborhood
$\Omega_\eta=\{x\in\Omega;\ \delta_\Omega(x)<\eta\}$
is $C^2$-smooth and 
$$\hbox{$\mu(x)\ge c_1\delta^\sigma_\Omega(x)$ on $\overline\Omega_\eta$.}$$
for some $\sigma<2$ and $c_1>0$. Without loss of generality,
we can assume $\sigma>1$.
By (\estimGraduDelta) we have
$$\int_{\Omega_\eta}\delta_{\Omega_\eta} ^{\sigma+1}\ |\nabla u|^2\le 
\int_{\Omega_\eta} \delta_\Omega^{\sigma+1} |\nabla u|^2\le C.
\leqno(\estimGradBB)$$
Applying Propostion 2(ii) with $p=2$, $a=1+\sigma$, $k=0$, 
this along with (\estimLoneMinus) guarantees
$$\int_{\Omega_\eta}  \delta_{\Omega_\eta} ^{\sigma-1}u^2 \le C.
$$
Next observe that
$\delta_{\Omega_\eta}(x)=\delta_\Omega(x)$ on $\Omega_{\eta/2}$, hence,
recalling $\sigma>1$,
$$\int_{\Omega_{\eta/2}}  \delta_{\Omega}^{\sigma-1}u^2 \le C.
\leqno(\estimHardyB)$$
On the other hand, by (\estimLoneMinus), we have, for all $\eps>0$,
$$\int_{\Omega_{\eta/2}}  \delta_{\Omega}^{\eps-1}u \le C(\eps).
\leqno(\estimLoneMinusB)$$

We now interpolate between (\estimHardyB) and (\estimLoneMinusB).
(It can be checked that working, instead of the Hardy-type estimate (\estimHardyB),
 with a weighted Sobolev estimate deduced from (\estimGradBB) and Propostion 2(i),
would not improve the conditions.)
To this end, we apply Lemma 3 
with $\omega=\Omega_{\eta/2}$, $\phi=\delta_\Omega$,
$$m={3\over 2}+\eps,\quad q=1,\quad r=2,\quad \gamma=0,
\quad b=1-\eps,\quad d=\sigma-1,$$
with $\eps>0$ small
(note that $q\le m \le r$). Taking $\eps\in (0,1/4)$ small enough, condition (\condGammaInterpol) is true due to
$$\sigma-1-\Bigl(2-{3\over 2}\Bigr)\sigma={\sigma\over 2}-1<0,$$
hence
$$\int_{\Omega_{\eta/2}} u^{(3/2)+\eps} \le C.$$
Since, owing to (\estimLpDelta), $n=3$ and $\eps<1/4$,
we also have
$$\int_{\Omega\setminus\Omega_{\eta/2}} u^{(3/2)+\eps} \le 
C\Bigl(\int_{\Omega\setminus\Omega_{\eta/2}} \delta_\Omega u^{2-\eps}\Bigr)^{((3/2)+\eps)/(2-\eps)} \le C,$$
we finally obtain
$$\int_\Omega c^{3/2}u^{{3\over 2}+\eps}\le C$$
and conclude by Proposition 1.
\fin

\medskip

{\bf 3. Proof of Theorem 5.} 
\medskip

For each integer $j\ge 1$, we seek $u=u_j$ under the form
$$u(x)=\cases{
jx,\quad 0\le x<1,\cr
\noalign{\vskip 1mm}
j+\log\bigl(1+j(x-1)\bigr),\quad 1\le x<2,\cr
\noalign{\vskip 1mm}
A_j\sin\bigl(\bigl({\pi\over 2}+\eps_j\bigr)(3-x)\bigr),\quad 2\le x\le 3,\cr
}
$$
where $A_j, \eps_j>0$ are to be determined.
The function $u$ satisfies 
$$-u''(x)=\cases{
0,\quad 0<x<1,\cr
\noalign{\vskip 1mm}
{u'}^2,\quad 1<x<2,\cr
\noalign{\vskip 1mm}
\bigl({\pi\over 2}+\eps_j\bigr)^2u,\quad 2<x<3.\cr
}
$$
Moreover $u(0)=u(3)=0$ and $u, u'$ are continuous at $x=1$.
It thus suffices to choose $A_j$ and $\eps_j$ in such a way as to ensure the continuity of $u$ and $u'$ at $x=2$.
We compute
$$u(2^-)=j+\log\bigl(j+1\bigr),\quad u(2^+)=A_j\sin\bigl(\textstyle{\pi\over 2}+\eps_j\bigr)=A_j\cos \eps_j$$
and
$$u'(2^-)={j\over j+1},\quad u'(2^+)=-\bigl(\textstyle{\pi\over 2}+\eps_j\bigr)A_j\cos\bigl(\textstyle{\pi\over 2}+\eps_j\bigr)
=\bigl(\textstyle{\pi\over 2}+\eps_j\bigr)A_j\sin\eps_j.$$
The conditions $u(2^-)=u(2^+)$ and $u'(2^-)=u'(2^+)$
are then equivalent to
$$A_j\cos \eps_j=j+\log\bigl(j+1\bigr).
\leqno(\EqnUnbddA)$$
and 
$$\bigl(\textstyle{\pi\over 2}+\eps_j\bigr)\tan \eps_j=\displaystyle{j\over j+1}\Bigl(j+\log\bigl(j+1\bigr)\Bigr)^{-1}
\leqno(\EqnUnbddB)$$
Since the function $s\mapsto \bigl(\textstyle{\pi\over 2}+s\bigr)\tan s$ is strictly increasing from $[0,\pi/2)$ onto $[0,+\infty)$,
there exists (a unique) $\eps_j\in (0,\pi/2)$ satisfying (\EqnUnbddB) and $A_j$ is then directly given by~(\EqnUnbddA).
Finally, we observe that $\|u_j\|_\infty\ge u_j(1)=j$ and that $\eps_j\to 0$ as $j\to\infty$.
The result is proved.
\fin

\medskip
\ \ {\bf Appendix: Proof of Proposition 2}
\medskip

By density (and Fatou's Lemma) it suffices to establish (\IneqSobolev) and (\IneqHardy)
for $v\in C^1(\overline\omega)$ with $v\ge 0$,
which we assume in the sequel.
In this proof, $C$ will denote a generic positive constant independent of $v$.
We proceed in several steps.

\medskip

{\bf Step 1.} Fix $1\le p<n$, $a>p-1$ and $k\ge 0$. Denote by $\lambda_1$ and $\varphi$ the first eigenvalue and eigenfunction of $-\Delta$
in $H^1_0(\omega)$. 
(If $\omega$ is not connected, we take $\varphi$ to be positive and $L^1$-normalized on each connected 
component of $\omega$.)
By elliptic regularity, we have $\varphi\in W^{2,q}(\omega)\cap C^2(\omega)$
for all $q\in (1,\infty)$. Then, as a consequence of the Hopf Lemma, there exist constants $c_1,c_2,\eta>0$,
such that
$$c_1\delta(x)\le \varphi(x)\le c_2\delta(x),\quad x\in\omega
\leqno(\equivHopf)$$
and
$$c_1\le |\nabla\varphi|\le c_2\quad\hbox{ on $\omega_\eta=\{x\in\omega;\ \delta(x)\le\eta\}$}.\leqno(\ControlNablaVarphi)$$
In particular, by (\equivHopf), it is equivalent to prove (\IneqSobolev) and (\IneqHardy)
 with $\varphi$ instead of $\delta$.

\medskip

{\bf Step 2.} We claim that
$$\int_\omega \delta^{a-p}v^{p}
\le C\int_\omega \delta^{a+2-p} v^p
+C\int_\omega \delta^a |\nabla v|^p.
\leqno(\ClaimAaa)$$

Set $\beta=a+1-p>0$. For $\eps>0$, we set $\psi=\varphi+\eps$
(the dependence on $\eps$ being omitted for brevity).
Using 
$$\Delta\Bigl({\psi^{1+\beta}\over 1+\beta}\Bigr)=
\nabla\cdot\Bigl(\psi^{\beta}\nabla\varphi\Bigr)
=\beta\psi^{\beta-1}|\nabla\varphi|^2+\psi^{\beta}\Delta\varphi
=\beta\psi^{\beta-1}|\nabla\varphi|^2-\lambda_1\psi^{\beta}\varphi$$
and noting that $v^p\psi^\beta\nabla\varphi\in (W^{1,q}(\omega))^n$
for any $q\in [1,\infty)$, we may apply the divergence theorem to deduce 
$$\eqalign{
I_\eps:=\eps^{\beta}\int_{\partial\omega}v^p{\partial\varphi\over\partial\nu}\,d\sigma &=\int_{\partial\omega}v^p\psi^{\beta}{\partial\varphi\over\partial\nu}\,d\sigma 
=\int_\omega \nabla\cdot\Bigl(v^p\nabla\Bigl({\psi^{1+\beta}\over 1+\beta}\Bigr)\Bigr)\,dx \cr
&=p\int_\omega v^{p-1}\psi^{\beta}\nabla v\cdot\nabla\varphi
+\int_\omega v^p\bigl[\beta\psi^{\beta-1}|\nabla\varphi|^2-\lambda_1\psi^{\beta}\varphi\bigr]. }$$
Therefore, by (\ControlNablaVarphi) and Young's inequality
(if $p>1$, or directly if $p=1$), we deduce
$$\eqalign{c_1^2\beta \int_{\omega_\eta} \psi^{\beta-1}v^p
&\le \beta \int_\omega \psi^{\beta-1}|\nabla\varphi|^2 v^p\cr
&\le I_\eps+\lambda_1  \int_\omega \psi^{\beta}\varphi v^p 
+p\int_\omega \bigl(\psi^{(\beta-1)(p-1)/p}v^{p-1}\bigr)\,
\bigl(\psi^{(\beta+p-1)/p}|\nabla v\cdot\nabla\varphi|\bigr)\cr
&\le I_\eps+\lambda_1  \int_\omega \psi^{\beta}\varphi v^p 
+{c_1^2\beta\over 2} \int_{\omega} \psi^{\beta-1}v^p
+C\int_\omega \psi^{\beta+p-1}|\nabla v|^p \cr
}$$
(here and below, $C$ is independent of $\eps$, as well as of $v$).
Consequently, 
$$
{c_1^2\beta\over 2} \int_{\omega_\eta} \psi^{\beta-1}v^p
\le I_\eps+\lambda_1  \int_\omega \psi^{\beta}\varphi v^p 
+{c_1^2\beta\over 2} \int_{\omega\setminus\omega_\eta} \psi^{\beta-1}v^p
+C\int_\omega \psi^{\beta+p-1}|\nabla v|^p.$$
Since $\int_{\omega\setminus\omega_\eta} \psi^{\beta-1}v^p
\le C \int_{\omega\setminus\omega_\eta} \psi^{\beta}\varphi v^p$
due to (\equivHopf), we deduce that
$$ \int_{\omega} (\varphi+\eps)^{\beta-1}v^p
\le I_\eps+C  \int_\omega (\varphi+\eps)^{\beta}\varphi v^p 
+C\int_\omega (\varphi+\eps)^{\beta+p-1}|\nabla v|^p.$$
Since $I_\eps\to 0$ as $\eps\to 0^+$, 
 (\ClaimAaa) follows by letting $\eps\to 0^+$ and using
 the monotone convergence theorem on the LHS 
 (which in particular yields the convergence of the integral on the LHS of (\ClaimAaa)).

\medskip

{\bf Step 3.} We next claim that, for any $r\ge 0$,
$$\int_\omega \varphi^{a-p}v^{p}
\le C\int_\omega \varphi^r v^p
+C\int_\omega \varphi^a |\nabla v|^p.
\leqno(\ClaimAab)$$

It suffices to prove (\ClaimAab) for $r=a-p+2j$ with $j$ integer $\ge 1$.
This is true for $j=1$ by Step 2.
Assume it is true for some integer $j\ge 1$.
By Step 2 (with $a$ replaced by $a+2j$), we have
$$\eqalign{\int_\omega \varphi^{a-p+2j}v^{p}
&\le C\int_\omega \varphi^{a-p+2(j+1)} v^p+C\int_\omega \varphi^{a+2j} |\nabla v|^p \cr
&\le C\int_\omega \varphi^{a-p+2(j+1)} v^p+C\int_\omega \varphi^{a} |\nabla v|^p.
}$$
Therefore,
$$\int_\omega \varphi^{a-p}v^{p}
\le C\int_\omega \varphi^{a-p+2j} v^p+C\int_\omega \varphi^a |\nabla v|^p
\le C\int_\omega \varphi^{a-p+2(j+1)} v^p+C\int_\omega \varphi^{a} |\nabla v|^p
$$
and (\ClaimAab) follows by induction.

\medskip

{\bf Step 4.} We claim that
$$
\Bigl(\int_\omega \varphi^{na/(n-p)} v^{p^*}\Bigr)^{p/p^*}
\le C\int_\omega\varphi^{a-p}v^{p}
+C\int_\omega\varphi^a|\nabla v|^{p}.
\leqno(\ClaimAa)$$

Let $s\ge 1$, $d>0$ to be chosen and let $\phi=\varphi^d$.
Note that $|\nabla\phi|=d|\nabla\varphi|\varphi^{d-1}\le C\varphi^{d-1}$.
It is well known that $\varphi^{d-1}\in L^1(\omega)$ (but this follows for instance
from (\ClaimAaa) with $p=1$, $a=d$ and $v\equiv 1$).
Therefore, $\phi v^{s}\in W^{1,1}(\omega)$.
Set $m=n'=n/(n-1)$. Observe that $p'>m$ due to $1<p<n$. By the standard Sobolev inequality in $W^{1,1}$, we have
$$
\int_\omega \phi^m v^{sm}
=\|\phi v^{s}\|_m^m\le C\|\nabla(\phi v^{s})\|_1^m+C\|\phi v^{s}\|_1^m,$$
hence
$$
\int_\omega \varphi^{dm} v^{sm}
\le C\|\phi v^{s-1}\nabla v\|_1^m+C\|v^{s}\nabla\phi\|_1^m+C\|\phi v^{s}\|_1^m.
\leqno(\ineqAa)$$
If $p=1$, then (\ClaimAa) follows from (\ineqAa) with the choice $s=1$, $d=a$.
Thus assume $p>1$.

By H\"older's inequality, for any $\theta\in [0,1]$, we have
$$\eqalign{
\|\phi v^{s-1}\nabla v\|_1^m
&\le C\Bigl(\int_\omega\phi^{(1-\theta)p'}v^{(s-1)p'}\Bigr)^{m/p'}
\Bigl(\int_\omega\phi^{\theta p}|\nabla v|^{p}\Bigr)^{m/p}.
}$$
Choosing 
$$s=(n-1)p/(n-p)>1,\quad
\theta=(n-p)/((n-1)p)\in (0,1),\quad
d=(n-1)a/(n-p),$$
 we get
$(1-\theta)p'=m$, $(s-1)p'=sm=p^*$, $\theta p=a/d$
and $(p'/m)'=p^*/m$. By Young's inequality, we then have
$$
\|\phi v^{s-1}\nabla v\|_1^m
\le {1\over 4}\int_\omega\varphi^{dm}v^{p^*}
+C\Bigl(\int_\omega\varphi^{a}|\nabla v|^{p}\Bigr)^{p^*/p}.
\leqno(\ineqAb)$$

Next, by H\"older's inequality, for any $q>1$, $\tau\in [0,1]$ and $\lambda\in\R$, we have
$$\eqalign{
\|v^{s}\nabla\phi\|_1^m+\|\phi v^{s}\|_1^m
&\le C\|\varphi^{d-1}v^{s}\|_1^m \cr
&\le C\Bigl(\int_\omega \varphi^{(d-1)\lambda q}v^{\tau sq}\Bigr)^{m/q}
\Bigl(\int_\omega \varphi^{(d-1)(1-\lambda)q'}v^{(1-\tau)sq'}\Bigr)^{m/q'}.
}
\leqno(\ineqAbb)$$
Choosing 
$$q=(ms-p)/(s-p)=p'>1,\quad
\tau=m/q=1-{n-p\over (n-1)p}\in (0,1),\quad
$$
 we get
$\tau sq=sm=p^*$ and
$(1-\tau)sq'=p$. If $d\neq 1$, we also choose
$\lambda=md/(d-1)q$
and we obtain
$(d-1)\lambda q=dm$ and
$$\eqalign{
(d-1)(1-\lambda)q'
&=(d-1)p\Bigl(1-{md\over (d-1)q}\Bigr) \cr
&=(d-1)p-md(p-1)={d(n-p)\over n-1}-p=a-p.
}$$
(If $d=1$, $\lambda$ is irrelevant in (\ineqAbb).)
 By Young's inequality, we then have
$$
C\|v^{s}\nabla\phi\|_1^m+C\|\phi v^{s}\|_1^m
\le {1\over 4}\int_\omega\varphi^{dm}v^{p^*}
+C\Bigl(\int_\omega\varphi^{a-p}v^p\Bigr)^{p^*/p}.
\leqno(\ineqAc)$$
Claim (\ClaimAa) now follows by combining (\ineqAa), (\ineqAb) and (\ineqAc)
(recalling $sm=p^*$ and $dm=na/(n-p)$).

\medskip

{\bf Step 5.} We claim that, for each $\ell\in [1,p]$, there exists
$b_\ell\ge 0$ (possibly depending also on $k$), such that
$$
\int_\omega \varphi^{b_\ell} v^\ell
\le C\Bigl(\int_\omega\varphi^kv\Bigr)^\ell
+C\Bigl(\int_\omega\varphi^a|\nabla v|^{p}\Bigr)^{\ell/p}.
\leqno(\ClaimBa)$$

When $\ell=1$, this is trivially true with $b_1=k$.

Next assume that (\ClaimBa) is true for some given $\ell\in [1,p]$.
It is then true for $\ell$ replaced with $\ell^*=n\ell/(n-\ell)$.
Indeed, letting $\bar a=\max(a,\ell+b_\ell)>\ell-1$, and
using (\ClaimAa) with $\ell, \bar a$ instead of $p, a$ and H\"older's inequality, we obtain
$$\eqalign{
\int_\omega \varphi^{n\bar a/(n-\ell)} v^{\ell^*}
&\le C\Bigl(\int_\omega\varphi^{\bar a-\ell}v^{\ell}\Bigr)^{\ell^*/\ell}
+C\Bigl(\int_\omega\varphi^{\bar a}|\nabla v|^{\ell}\Bigr)^{\ell^*/\ell}\cr
&\le C\Bigl(\int_\omega\varphi^{b_\ell}v^{\ell}\Bigr)^{\ell^*/\ell}
+C\Bigl(\int_\omega\varphi^a|\nabla v|^{p}\Bigr)^{\ell^*/p}\cr
&\le C\Bigl(\int_\omega\varphi^kv\Bigr)^{\ell^*}
+C\Bigl(\int_\omega\varphi^a|\nabla v|^{p}\Bigr)^{\ell^*/p}
}$$
that is, (\ClaimBa) with $b_{\ell^*}=n\bar a/(n-\ell)$.
By H\"older's inequality, Property (\ClaimBa) is then true
for any $\ell\in[1,\ell^*]$.

Now, for any integer $i\in\{0,\cdots,n-1\}$, set $p_i=n/(n-i)$,
and observe that $p_{i+1}:=p_i^*=np_i/(n-p_i)$ for $i\in\{0,\cdots,n-2\}$
(due to $1/p_i=1-(i/n)$).
Next let $j\in\{0,\cdots,n-2\}$ be the unique integer such that
$p_j\le p<p_{j+1}$.
By induction, starting from $\ell=p_0=1$,
 the previous paragraph guarantees that (\ClaimBa) is true for
any $\ell\in [1,p_{j+1}]\supset [1,p]$ and the claim follows.

\medskip
{\bf Conclusion.}
Combining (\ClaimAab) for $r=b_p$ and (\ClaimBa) for $\ell=p$, we obtain
$$
\int_\omega \varphi^{a-p}v^{p}
\le C\int_\omega \varphi^{b_p} v^p+C\int_\omega \varphi^a |\nabla v|^p
\le C\Bigl(\int_\omega \varphi^k u\Bigr)^p+C\int_\omega \varphi^a |\nabla v|^p,
$$
which proves assertion (ii).
Assertion (i) then follows from (\ClaimAa) and assertion (ii).
\fin

\vskip 6mm
{\baselineskip=11pt \parindent=1cm
\font\pc=cmcsc9
\font\rmn=cmr9
\font\sln=cmsl9
\font\rmb=cmbx8 scaled 1125 \rm
\font\it=cmti9

\rmn \eightpoint

\bigskip
\centerline{\bf REFERENCES}
\medskip

\item{[\ADP]} 
B. Abdellaoui, A. DallÕAglio, I. Peral, 
Some remarks on elliptic problems with critical growth in the gradient, 
{\sln J. Differential Equations} 222 (2006), 21-62 $\&$ Corr. 
{\sln J. Differential Equations} 246 (2009), 2988-2990.

\smallskip

\item{[\AP]} N. Alaa, M. Pierre, Weak solutions of some quasilinear elliptic equations with data measures, 
{\sln SIAM J. Math. Anal. 24}  (1993), 23-35.

\smallskip

\item{[\AL]} A. Alvino, P.L. Lions, G. Trombetti, Comparison results for elliptic and parabolic equations via Schwarz symmetrization, 
{\sln Ann. Inst. H. Poincar\'e, Analyse non lin\'eaire}  7(1990), 37-65.

\smallskip

\item{[\ACJT]} D. Arcoya, C. de Coster, L. Jeanjean, K. Tanaka,
Continuum of solutions for an elliptic problem with critical growth in the gradient,
preprint ArXiV 1304.3066 (2013).

\item{[\Av]} A. Avantaggiati,
On compact embedding theorems in weighted Sobolev spaces,
{\sln Czech. Math. J.} 29 (1979), 635-648.

\smallskip

\item{[\BVV]} M.-F. Bidaut-V\'eron, L. Vivier,
An elliptic semilinear equation with source term
involving boundary measures: the subcritical case
{\sln Rev. Mat. Ibero\-americana}
16 (2000), 477-513.

\smallskip

\item{[\BMPa]} L. Boccardo, F. Murat, J.-P. Puel, Existence of bounded solutions for nonlinear elliptic unilateral problems, 
{\sln Ann. Mat. Pura Appl.}  152 (1988), 183-196.

\smallskip

\item{[\BMPb]} L. Boccardo, F. Murat, J.-P. Puel, $L^\infty$ estimate for some nonlinear elliptic partial differential equations and application to an existence result, 
{\sln SIAM J. Math. Anal.}  23 (1992), 326-333.

\smallskip

\item{[\BC]} H. Brezis, X. Cabr\'e,
Some simple nonlinear PDEÕs without solutions, 
{\sln Boll. Unione Mat. Ital.}  (8) 1-B (1999), 223-262.

\smallskip

\item{[\BT]} H. Brezis, R.E.L. Turner, On a class of superlinear elliptic problems, 
{\sln Comm. Partial Differ. Equations}  2 (1977), 601-614.

\smallskip

\item{[\CRT]} M.G. Crandall, P.H. Rabinowitz, L. Tartar, On a Dirichlet problem with a singular nonlinearity, 
{\sln Comm. Partial Differential Equations}  2 (1977), 193-222.

\smallskip

\item{[\FM]} V. Ferone, F. Murat, Nonlinear problems having quadratic growth in the gradient: an existence result when the source term is small, 
{\sln Nonlinear Anal. TMA}  42 (2000), 1309-1326.

\smallskip

\item{[\FPa]} V. Ferone, M.R. Posteraro, On a class of quasilinear elliptic equations with quadratic growth in the gradient, 
{\sln Nonlinear Anal. TMA}  20 (1993), 703-711.

\smallskip

\item{[\FPb]} V. Ferone, M.R. Posteraro, J.M. Rakotoson, $L^\infty$-estimates for nonlinear elliptic problems with p-growth in the gradient, 
{\sln J. Ineq. Appl.}  2 (1999), 109-125.

\smallskip

\item{[\FSW]} M. Fila, Ph. Souplet, F. Weissler,
Linear and nonlinear heat equations in
$L^q_\delta$ spaces and universal bounds for global solutions, 
{\sln Math. Ann.}  320 (2001), 87-113.

\smallskip

\item{[\FMT]} S. Filippas, V. Maz'ya, A. Tertikas,
Critical Hardy-Sobolev inequalities. 
{\sln J. Math. Pures Appl.}  (9) 87 (2007), 37-56.

\smallskip

\item{[\GM]} N. Grenon-Isselkou, J. Mossino, Existence de solutions born\'ees pour certaines \'equa\-tions elliptiques quasilin\'eaires, 
{\sln C. R. Math. Acad. Sci. Paris } 321 (1995), 51-56.

\smallskip

\item{[\JS]} L. Jeanjean, B. Sirakov, Existence and multiplicity for elliptic problems with qua\-dratic growth in the gradient, 
{\sln Comm. Part. Diff. Equ.} 38 (2013), 244-264.

\smallskip

\item{[\MPS]} C. Maderna, C. Pagani, S. Salsa, Quasilinear elliptic equations with quadratic growth in the gradient, 
{\sln J. Differential Equations}  97 (1992), 54-70.

\smallskip

\item{[\Ne]} J. Ne\v cas,
Sur une m\'ethode pour r\'esoudre les \'equations aux d\'eriv\'ees partielles du type elliptique, voisine de la variationnelle,
{\sln Ann. Scuola Normale Sup. Pisa} 16 (1962), 305-326.

\smallskip

\item{[\OpKu]} B. Opic, A. Kufner, 
Hardy-type inequalities. 
Pitman Research Notes in Mathematics Series, 219. Longman Scientific \& Technical, Harlow, 1990. 

\smallskip

\item{[\QSa]} P. Quittner, Ph. Souplet,
A~priori estimates and existence for elliptic systems via
bootstrap in weighted Lebesgue spaces,
{\sln Arch. Rational Mech. Anal.} 174 (2004), 49-81.

\smallskip

\item{[\QSb]} P. Quittner, Ph. Souplet,
Superlinear parabolic problems. Blow-up, global existence and steady states,
Birkhauser Advanced Texts, 2007, 584 p.+xi. ISBN: 978-3-7643-8441-8

\bye